\documentclass
[12pt,a4paper]{article}
\usepackage[cp1251]{inputenc}
\usepackage{amssymb, amsfonts, amsmath, latexsym}
\usepackage[english]{babel}

\begin{document}

\normalsize\centerline{\bf  I.V. Protasov and K.D. Protasova}\vspace{6 mm}

\Large
\centerline{\bf  Lattices of coarse structures }\vspace{6 mm}

\normalsize


{\bf Abstract.} We consider the lattice of coarse structures on a set $X$ and  study metrizable, locally finite and  cellular coarse structures on $X$ from the lattice  point of view.

\vspace{6 mm}

UDC 519.51

\vspace{3 mm}

\vspace{6 mm}

\section{Introduction}

Following [11], we say that a family $\mathcal{E}$ of subsets of $X\times X$ is a {\it coarse structure} on a set $X$ if

\begin{itemize}
\item{} each $\varepsilon \in \mathcal{E}$ contains the diagonal $\vartriangle _{X}$, $\vartriangle _{X}= \{(x,x) : x \in X\}$ ; \vskip 5pt

\item{}  if $\varepsilon, \delta\in\mathcal{E}$ then $\varepsilon \circ\delta\in\mathcal{E}$  and $\varepsilon^{-1}\in\mathcal{E}$ where $\varepsilon \circ\delta = \{(x, y): \exists z ((x,z)\in\varepsilon, (z,y)\in\delta)\}, $  $ \ \varepsilon^{-1}= \{(y,x): (x,y)\in\varepsilon\}$;

\item{}  if $\varepsilon\in\mathcal{E}$ and $\bigtriangleup_{X}\subseteq \varepsilon^{\prime}\subseteq\varepsilon$ then $\varepsilon^{\prime}\in\mathcal{E}$.

\end{itemize}

Each $\varepsilon\in\mathcal{E}$ is called an {\it entourage} of the diagonal. We note that
$\mathcal{E}$ is closed under finite union
$(\varepsilon\subseteq\varepsilon\delta, \ \delta\subseteq \varepsilon\delta)$,  but  $\mathcal{E}$ is not an ideal in the Boolean algebra of all subsets
of $X\times X$  because $\mathcal{E}$ is not closed under formation of all subsets of its members.

A subset $\mathcal{E}^{\prime}\subseteq\mathcal{E}$ is called {\it a base} for $\mathcal{E}$ if, for
every $\varepsilon\in\mathcal{E}$ there exists $\varepsilon^{\prime}\in\mathcal{E}^{\prime}$ such that $\varepsilon\subseteq\varepsilon^{\prime}$.

The pair $(X, \mathcal{E})$ is called a {\it coarse space}. For
$x\in X$  and $\varepsilon\in\mathcal{E}$, we denote $B(x, \varepsilon)= \{y\in X: (x,y)\in\varepsilon \}$
and say that  $B(x,\varepsilon)$ is a {\it ball of radius  $\varepsilon$ around $x$.}
 We note that a coarse space can be considered as
 an asymptotic counterpart of a uniform topological space and could be defined in terms of balls, see
 [7], [9]. In this case a coarse space is called a {\it ballean}.
  For categorical look at the  balleans and coarse structures as two faces of the same
  coin  see [1].

A coarse structure $\mathcal{E}$ on $X$ is called {\it connected} if, for any $x, y \in X$, there is $\varepsilon\in\mathcal{E}$ such that $y\in B(x,\varepsilon)$.
 A subset $Y$ of $X$ is called {\it bounded} if there exist $x\in X$ and $\varepsilon\in\mathcal{E}$
   such that $Y\subseteq B(x, \varepsilon)$.  A coarse
   structure $\mathcal{E}$ is called {\it bounded} if $X$ is bounded,
   otherwise $\mathcal{E}$ is called    {\it unbounded}. We note that
    on every set $X$,  there exists the unique connected bounded coarse structure $\{\varepsilon\subseteq X\times X: \bigtriangleup_{X}\subseteq\varepsilon\}$, and
    if $X$ is finite this structure is the unique connected coarse  structure on $X$.

In what follows, we consider only {\bf connected}  coarse structures on {\bf infinite} sets.

Given a set $X$, the family $\mathcal{L}_{X}$  of all coarse  structures
on $X$ is partially ordered by the inclusion,
and $\mathcal{L}_{X}$  can be considered as a lattice  with the operations $\wedge$ and  $\vee$:

$\bullet$  $\mathcal{E}\wedge \mathcal{E}^{\prime}$  is the strongest coarse  structure such that
$\mathcal{E}\wedge \mathcal{E}^{\prime}\subseteq\mathcal{E}$, $\mathcal{E}\wedge \mathcal{E}^{\prime}\subseteq\mathcal{E}^{\prime}$;

$\bullet$  $\mathcal{E}\vee \mathcal{E}^{\prime}$  is the weakest coarse structure such that
$\mathcal{E}\subseteq\mathcal{E}\vee \mathcal{E}^{\prime}$, $\mathcal{E}^{\prime}\subseteq\mathcal{E}\vee \mathcal{E}^{\prime}$.

More explicitly, we have

$(\wedge) \   \   \  \ $   $   \   \  \mathcal{E}\wedge \mathcal{E}^{\prime}=\{\varepsilon\cap  \varepsilon^{\prime}:  \varepsilon\in\mathcal{E},  \  \varepsilon^{\prime}\in\mathcal{E}^{\prime}\} . $

To clarify  the operation $\vee$, we note that
$\mathcal{E}\vee \mathcal{E}^{\prime}$
is the smallest coarse structure on $X$ coutaing
$\mathcal{E}$ , $\mathcal{E}^{\prime}$  and closed under the operations
$\circ$ and $^{-1}$ and  taking subsets of its members containing  $\triangle_{X}$.
Since $(\varepsilon  \circ \delta)^{-1}=\delta ^{-1}\circ\varepsilon^{-1}$, the usage of the inversion is superfluous. Hence, we have

$(\vee)  \  \  \  \  $  $  \  \mathcal{E} \vee  \mathcal{E}^{\prime}$
has the base  $\{\varepsilon_{1}\circ\ldots\circ\varepsilon_{n}:  n\in\mathbb{N}, \  \varepsilon_{i}\in\mathcal{E}\bigcup\mathcal{E}^{\prime}\}$.

The coarse structure ${\bf 1} _{X}=\{\varepsilon: \triangle_{X}\subseteq\varepsilon\}$
is the unit of
the lattice $\mathcal{L_{X}}$, and the coarse structure
 ${\bf 0} _{X}=\{\varepsilon: \triangle_{X}\subseteq\varepsilon\}$, $\varepsilon\setminus\bigtriangleup_{X}$  is  finite $\}$ is the null of $\mathcal{L}_{X}$.

The purpose of this note is to explore lattice properties of some basic coarse structures on a set $X$, namely, metrizable, locally finite and cellular structures. The main results are exposed in section 4 and concern $\vee$-decomposability of coarse structure.


\section{Metrizable structures}

Every metric $d$ on a set $X$  defines the coarse
structure $\mathcal{E}_{d}$ on $X$ with the base  $\{(x,y): d(x,y)\leq n\}$,  $n\in \omega$.
A coarse structure $\mathcal{E}$ on $X$ is called {\it metrizable} if
there exists a metric $d$ on $X$ such  that $\mathcal{E}= \mathcal{E}_{d}$.
By [9, Theorem 2.1.1], $\mathcal{E}$ is metrizable if and only if
$\mathcal{E}$ has a countable base. In view of $(\wedge)$ and $(\vee)$,
the family $\mathfrak{M} _{X}$  of all metrizable structures   on $X$ is a
sublattice of $\mathcal{L} _{X}$.
Since ${\bf 1}_{X}$ is metrizable, $\mathfrak{M} _{X}$  has
the unit. If $X$ is countable then
${\bf O} _{X}$ is metrizable so
$\mathfrak{M} _{X}$ has the null.
In Proposition 1, we  show that $\mathfrak{M}$  does not have null if $X$ is uncountable.

We use the following simple construction. Let
$(X,d)$ and  $(Y,\rho)$    be metric spaces such that $X\cap Y=\emptyset$.
We take  $x_{0}\in X$, $y_{0}\in Y$,  $r>0$ and define a
metric $\mu$ on $X\bigcup Y$ by the following rule:
  $\mu (x, x^{\prime})= d(x, x^{\prime})$
  if $x, x^{\prime} \in X$; $\mu(y, y^{\prime})=\rho (y, y^{\prime})$
  if $(y, y^{\prime})\in Y$;
  $\mu (x, y)= d(x, x_{0}) + r + \rho(y, y_{0})$
  if $x\in X$, $y\in Y$.
The obtained metric space  $(X\bigcup Y,\mu)$ is called the
{\it $(x_{0}, y_{0}, r)$-join} of $(X,d)$  and $(Y,\rho)$.
\vspace{6 mm}

{\bf Proposition 1. } {\it If  an unbounded metric space $(X, d)$ contains an infinite bounded subset $Y$ then there exists a metric $\mu$  on $X$ such that $\mathcal{E}_{\mu}\subset \mathcal{E} _{d}$.
\vspace{6 mm}

Proof.} We pick $m\in \mathbb{N}$ such that $d(y, y^{\prime})<m$ for
all $y,y^{\prime} \in Y$ and choose an unbounded metric $\rho$
such that   $\rho(y,y^{\prime})>m$
 for all $y,y^{\prime} \in Y$ .
 Then we take $x_{0}\in X$,  $y_{0}\in Y$, put $r=d(x_{0}, y_{0})$ and consider the
  $(x_{0}, y_{0}, r)$-join $(X,\mu)$ of
  $(X\setminus Y, d)$ and $(Y,\rho)$. By
  the choice of $\mu$, $\mathcal{E}_{\mu}\subseteq \mathcal{E}_{d}$. Since $Y$ is unbounded in
  $(X,\mu)$, we have $\mathcal{E}_{\mu}\subset \mathcal{E}_{d}$. $ \  \Box$
\vspace{6 mm}

If $X$ is uncountable then some ball in $(X,d)$ must be infinite and, applying Proposition 1, we get
$\mathcal{E}_{\mu}$ such that $\mathcal{E}_{\mu}\subset \mathcal{E}_{d}$.
\vspace{6 mm}

{\bf Proposition 2. } {\it
For every unbounded metric space $(X, d)$, there exists an unbounded metric $\rho$ on $X$ such  that  $\mathcal{E}_{d}\subset \mathcal{E}_{\rho}$.
\vspace{6 mm}

Proof. } We choose a sequence $(a_{n})_{n\in\omega} $  in $X$ such that
$B_{d}(a_{n},n)\cap B_{d}(a_{m},n)=\emptyset$ for  all $n>m$.
Then we denote $Y=\{a_{2n}: n\in\omega\}$, $\delta=\bigtriangleup_{X}\cup\{(x,y): x,y \in Y\}$
 and put
 $$\mathcal{H}= \{\varepsilon_{0}\circ \ldots\circ\varepsilon_{n}: \varepsilon_{i}\in \mathcal{E}_{d}\cup\{\delta\}\}.$$
Then $\mathcal{H} $ is the base for some uniquely defined metrizable coarse structure $\mathcal{E}$ on $X$ so $\mathcal{E}= \mathcal{E}_{\rho}$ for some metric $\rho$ on $X$.

We note that $\mathcal{E}_{d}\subset \mathcal{E}$
  because $Y$ is bounded in $\mathcal{E}$ but $Y$ is unbounded in $\mathcal{E}_{d}$.

To prove that $\mathcal{E}$ is unbounded, we choose $y_{0}\in Y$ and
observe by induction on $n$ that, for each ball
$B(y_{0}, \varepsilon_{0}\circ\ldots\circ\varepsilon_{n})$, $\varepsilon_{i}\in\mathcal{E}$,  there exists
   $m\in \omega$ such that
$B(y_{0}, \varepsilon_{0}\circ\ldots\circ\varepsilon_{n})\subseteq B_{d}(Y, m)$   where $B_{d}(Y, m)= \cup_{y\in Y} B_{d}(y,m)$.
By the choice of $(a_{n})_{n\in\omega}$, $a_{k}\notin B_{d}(Y,m)$ for each odd $k>m$. $ \  \Box$
\vspace{6 mm}

{\bf Proposition 3. } {\it
Let $d,\mu$  be unbounded metrics on a set $X$ such that $\mathcal{E} _{d}\subset \mathcal{E}_{\mu}$. Then there exists a metric $\rho$ on $X$ such that $\mathcal{E} _{d}\subset\mathcal{E}_{\rho}\subset \mathcal{E}_{\mu}$.

  Proof.}  Since  $\mathcal{E}_{d}\subset \mathcal{E}_{\mu}$, we can choose $m\in\omega$  and
   $a_{n},  c_{n}, n\in\omega$  such that

   $(1)   \  \   d(a_{n}, c_{n})>n;$

    $(2)   \  \   \mu(a_{n}, c_{n})< m.$
\vspace{4 mm}

We denote $A=\{a_{n} :  n \in\omega\}$,  $C=\{c_{n} :  n \in\omega\}$ and consider two cases.
\vspace{4 mm}

{\it Case: } $A\cup C$ is unbounded in $(X, \mu)$. We assume that $C$ is unbounded. Passing to subsequences, we may suppose  that

$(3)   \  \   B_{\mu}( c_{n}, n )\cap B_{\mu}( c_{k}, k )=\emptyset$  for all distinct $n, k \in \omega$.

We put $\delta=\bigtriangleup_{X} \cup \{(a_{2n}, c_{2n}),  (c_{2n},  a_{n}): n\in\omega\}$ and a
  coarse  structure $\mathcal{E}$ on $X$  with the base
  $$\mathcal{H}= \{\varepsilon_{0}\circ \ldots\circ\varepsilon_{n}: \varepsilon_{i}\in \mathcal{E}_{d}\cup\{\delta\}\}.$$

Clearly, $\mathcal{E}$ is metrizable so  $\mathcal{E} = \mathcal{E}_{\rho}$   for some metric
$\rho$ on $X$. Since $\mathcal{E} _{d}  \subseteq \mathcal{H}$,  we  have $\mathcal{E}_{d} \subseteq \mathcal{E}_{\rho}$.  By (1), $\mathcal{E}_{d} \subset \mathcal{E}_{\rho}$.

By the choice of $\delta$ and (2), $\mathcal{E}_{\rho} \subseteq \mathcal{E}_{\mu}$.
We take an arbitrary $\varepsilon_{0}\circ \ldots\circ\varepsilon_{n}\in\mathcal{H}$,
 delete all $\varepsilon_{i},  \  \varepsilon_{i}=\delta$
  and get  $\varepsilon_{0}^{\prime}\circ \ldots\circ\varepsilon_{s}^{\prime}$.
  Since $\varepsilon_{0}^{\prime}\in \mathcal{E}_{d},\ldots, \varepsilon_{s}^{\prime}\in \mathcal{E}_{d}$,
  there exists
  $t\in \omega$ such that
  $B  (x, \varepsilon_{0}^{\prime}\circ \ldots\circ\varepsilon_{s}^{\prime})\subseteq  B_{d}(x,t)$
   for each $x\in X$.  If $2k+1>t$  then by (1), (3)  and the choice of $d$, we have
   $a_{2k+1}\notin B(c_{2k+1}, \varepsilon_{0}\circ \ldots\circ\varepsilon_{n})$.
     By (2), $a_{2k+1}\in B_{\mu}(c_{2k+1}, m).$
Hence, $\mathcal{E}_{d} \subset  \mathcal{E}$.
\vspace{4 mm}

{\it Case:} $A\cup B$  is bounded $(X, \mu )$.  By (1),  either $A$ or $C$  is unbounded  in  $(X, d)$.
We assume that $A$ is  unbounded  in  $(X, d)$.
Passing to subsequences, we may suppose
that $B_{d}(a_{n}, n)\cap  B_{d}(a_{k}, n)=\emptyset$
for all  $n> k$.  Since $A$ is bounded in $(X, \mu)$, we can apply
arguments from the proof of   Proposition 2. $ \  \Box$


\section{Locally finite structures }

A coarse structure  $\mathcal{E}$  on a set $X$ is called
 {\it locally finite} if each ball $B(x,\varepsilon)$,  $x\in X$,  $\varepsilon\in \mathcal{E}$ is finite.
 If  for every  $\varepsilon\in \mathcal{E}$
 there exists $m\in \omega$ such that,
 for each $x\in X$,  $|B(x, \varepsilon )|< m$  then $\mathcal{E}$ is called
{\it uniformly } locally finite.
The uniformly locally  finite coarse structures play the central part in {\it  Asymptotic Topology}, see [2], [6].

In view $(\wedge)$ and $(\wedge)$,  the families of  locally finite and uniformly locally finite  coarse structures on $X$ are sublattices   of $\mathcal{L} _{X}$  with the null ${\bf 0}_{X}$. To find units in
 these structures, we use some specializations  of the following general construction.

Let $G$ be a group with the identity $e$ and let $X$ be a   $G$-space with the
action $G \times X \longrightarrow X$,    $ \  (g,x) \longmapsto gx$. An ideal $\mathcal{I}$  in the Boolean
algebra $\mathcal{P}_{G}$  of all subsets of $X$ is called a {\it  group ideal}
 if $\mathcal{I}$ contains the ideal $\mathfrak{F}_{G}$ of all finite subsets and
 $AB ^{-1} \in \mathcal{I}$  for any $A, B\in \mathcal{I}$. For group ideals see  [8] and  [9,  Chapter 6].

We assume that $G$ acts on $X$ transitively  (for any $x,y\in X$,  there is $g\in G$  such that $gx=y$) and, given a group  ideal $\mathcal{I}$  in  $\mathcal{P}_{G}$ ,  we consider the coarse structure $\mathcal{E},  \mathcal{E}= \mathcal{E}(G, \mathcal{I}, X)$  on $X$  with the base   $\{\varepsilon_{A}: A \in \mathcal{I},  c\in A\}$,  $\varepsilon_{A} =\{(x, gx): x\in X, g\in A \}$. Then
$B (x,\varepsilon_{A})= Ax$, $Ax =\{gx: g\in A \}$.

By [5, Theorem 1],  for every coarse structure $\mathcal{E}$ on $X$ there exists a group $G$ of permutations of $X$  and a group ideal $\mathcal{I}$ in $\mathcal{P}_{G}$ such that $\mathcal{E}=\mathcal{E}(G, \mathcal{I}, X)$.

Every coarse structure  $\mathcal{E}(G, \mathfrak{F}_{G}, X)$  is uniformly locally finite.
 If $X=G$ and $G$  acts of $X$ by the  left shifts then  $\mathcal{E}(G, \mathfrak{F}_{G}, G)$ is called the  {\it finitary  coarse structure on} $G$.
 The finitary coarse structures on finitely generated groups take an impotent part of {\it Geometrical Group Theory},  see [4, Chapter 4].

By  [6, Theorem 1],  for every uniformly locally finite coarse structure  $\mathcal{E}$ on $X$,  there exists a group $G$ of permutations of $X$ such that $\mathcal{E}= \mathcal{E}(G, \mathfrak{F}_{G}, X)$.
Applying this theorem, we get
\vspace{6 mm}

{\bf Proposition 4.} {\it Let $X$ be a set and let $S_{X}$ denotes the group of all  permutations of $X$.
Then the coarse structure $\mathcal{E}(S_{X}, \mathfrak{F}_{S_{X}}, X)$
 is the unit in the lattice of all uniformly locally finite coarse structures on $X$.}
 \vspace{6 mm}

 We endow $S_{X}$ with the topology of pointwise convergence, the subsets
 $\{g\in S_{X} : gx = y\}$,    $x,y\in X$ form a  sub-base of this topology.
  We note that a subset  $A$  of $S_{X}$
   is precompact if and only if the set $Ax$ is finite for each $x\in X$. We denote by $ \mathfrak{C} _{S_{X}}$ the ideal of all precompact   subsets of  $S_{X}$.

\vspace{6 mm}

{\bf Proposition 5.} {\it The coarse structure $\mathcal{E}(S_{X}, \mathfrak{C} _{S_{X}}, X)$
is
the unit in the latiice of all locally finite coarse structures on $X$.
\vspace{6 mm}

Proof.} Let ${\mathcal{E}}^{\prime}$  be  a locally finite and let $\varepsilon\in {\mathcal{E}}^{\prime}$.
 For
each $(x,y)\in \mathcal{E}$, we take a permutation
$f _{(x,y)}$ such that $f(x)=y$, $f(y)=x$ and
$f(z)=z$
for all  $z\in X \setminus \{x,y\} $.
We put $A=\{f_{x,y} : (x,y)\in\varepsilon\}$.
Since  $\mathcal{E}$  is  locally finite, we have
$A\in \mathfrak{C} _{S_{X}}$.
 By the choice of $A$,
 $\varepsilon=\varepsilon_{A}$ so $\varepsilon\in\mathcal{E}(S_{X}, \mathfrak{C} _{S_{X}}, X)$. $ \  \Box$
\vspace{4 mm}

We recall that two metric space $(X,d), (Y,\rho )$ is {\it isometric}
if there is a bijection $f: X\longrightarrow  Y$  such that
$d(x, y)= \rho (f(x), f(y))$
 for all  $x, y \in Y$. If a metrizable
 coarse structure on $X$ is locally finite  then $X$ must be countable.

\vspace{6 mm}

{\bf Proposition 6.} {\it
Let $d,\rho$  be unbounded metrics on $X$ such
 that $(X, d), (X,\rho)$ are locally finite. Then there
  exists a metric $\mu$  on $X$ such  that
   $\mathcal{E} _{d} \wedge \mathcal{E}_{\mu}  =  {\bf 0} _{X}$
and $(X,\rho),  (X, \mu)$
 are isometric.
\vspace{4 mm}

Proof.} Assume that we have defined $\mu$ so that,
for every $n\in\omega$ there exists a finite subset  $F$ of $X$
such that if $x\in X\setminus  F$, $y\in X, y\neq x$
 and
$d(x,y)<n$  then  $\mu(x,y)> n$.
Then
$B_{d}(x,n)\cap  B_{\mu}(x,n)= \{x\}$
and   $\mathcal{E} _{d} \wedge \mathcal{E}_{\mu}  =  {\bf 0} _{X}$
The problem is to find $\mu$  so  that $(X,\rho),(X, \mu)$
are isometric and $\mu$  satisfies above conditions.
 We construct inductively some  bijection  $f: X \longrightarrow X$ and  get $\mu$  as
 $\mu(x,y)=\rho(f(x), f(y) )$.

We choose a sequence $(a_{n})_{n\in\omega}$ in $X$  such that
$B_{d}(a_{n},n)\cap  B_{d}(a_{m},m)= \emptyset$  for all distinct $n,m\in \omega$.
Then we take a numeration $X= \{ x_{n} : n \in \omega\}$ satisfying
the following conditions.
If $x_{i}\in B_{d}(a_{n},n)$, $x_{j}\in B_{d}(a_{m},m)$
and $n< m$ then  $i< j$. For every  $n\in \omega$,
$B_{d}(a_{n},n)= \{ x_{k_{n}}, x_{k_{n}+1},  \ldots , x_{s_{n}} \}$, $x_{k_{n}}=a_{n}$.

We choose a sequence $(b_{n})_{n\in\omega}$ in $X$  such that
$B_{\rho}(b_{n},n)\cap  B_{\rho}(b_{m},m)= \emptyset$  for all distinct $n,m\in \omega$.

We put $f(x_{0})=b_{0}$
and suppose that we have defined  $f(x_{0}), \ldots ,f(x_{k})$. If $x_{k+1}\in X\setminus \{a_{n}: n\in\omega\}$ then we put  $f(x_{k+1})=b_{k+1}$.

If  $x_{k}=a_{n}$ and
$X\setminus (\{f(x_{0}),\ldots , f(x_{k})\}\cup  \cup_{i>k} B_{\rho}(b_{i}, i) )=\emptyset$
then we put $f(x_{k+1})=b_{k+1}$.
Otherwise, we take the minimal $j$ such that
$x_{j}\in X\setminus (\{f(x_{0}),\ldots , f(x_{k})\}\cup  \cup_{i>k} B_{\rho}(b_{i}, i) )$
  and put $f(x_{k+1})=x_{j}$.

After $\omega$ steps, we get the injective mapping $f: X\longrightarrow X$.
 To see that $f$ is surjective, it suffices note that
$\cap_{k\in\omega}\cup_{i>k}B_{\rho}(b_{i},i)=\emptyset$.

At last, if $n\in\omega$, $k>n$, $y\neq x_{k}$  and $d(x_{k}, y)<n$, $y\neq x_{k}$ then
$f(y)\notin B_{\rho}(f(x_{k}), n)$
 so $\rho (f(x_{k}), f(y))>n$. $ \  \Box$

\vspace{4 mm}

{\bf Proposition 7.} {\it Let $(X, d)$  be a locally finite metric space. Then there exists a metric $d^{\prime}$ on $X$  such that $(X, d^{\prime})$ is locally finite and $E_{d} \subset  E_{d^{\prime}}$.

\vspace{4 mm}

Proof.} If $E_{d}= {\bf 0 } _{X}$  then the statement is evident.
 We suppose that $E_{d}\neq  {\bf 0 } _{X}$, put $\rho=d$ and applying  Proposition 6 get
  the locally finite metric $\mu$ on  $X$  such that $\mathcal{E}_{\mu}\wedge  \mathcal{E}_{d}={\bf 0 } _{X}$
   and $(X, \mu)$ is isometric to $(X, d)$. Since $\mathcal{E}_{d} \neq {\bf 0 } _{X}$, we have
   $\mathcal{E}_{d} \subset \mathcal{E}_{d}\vee \mathcal{E}_{\mu}$  and
   $ \mathcal{E}_{d}\vee \mathcal{E}_{\mu}= \mathcal{E}_{d^{\prime}} $
    for some  locally finite metric $d^{\prime}$ on  $X$.   $ \  \Box$

\vspace{4 mm}

{\bf Proposition 8.} {\it Let  $(X, d), (X,\mu)$  be locally finite
metric spaces such that
$\mathcal{E}_{d}\subset \mathcal{E}_{\mu}$.  Then
there  exists a locally finite metric $\rho$ on $X$ such that
$\mathcal{E}_{d}\subset \mathcal{E}_{\rho} \subset \mathcal{E}_{\mu}$.

\vspace{4 mm}

Proof.} We apply Proposition 3 and note that $\rho$ is locally finite because
$ \mathcal{E}_{\rho} \subset \mathcal{E}_{\mu}$.
$ \  \Box$

\vspace{4 mm}

By Proposition 7, the unit in the lattice of all locally  finite coarse structures on $X$
(see  Proposition 5) is not metrizable. It is easy to see that Proposition  7
remains true with uniformly  locally finite metrics in  place of locally finite metrics,
 so the unit in the lattice of all uniformly locally finite
  metric spaces (see  Proposition 4) is not metrizable.


\section{ Cellular structures}

Every equivalence $\varepsilon$  on a set  $X$  partitions $X$ into classes of equivalence. On the other  hand, every partition  $\mathcal{P}$  of $X$  defines  the equivalence $\varepsilon$ by the rule: $(x,y)\in \varepsilon$ if
and only if $x,y$ are in some cell of the partition $\mathcal{P}$.

A coarse structure $\mathcal{E}$ on a set  $X$  is called {\it cellular}  if $\mathcal{E}$  has a base consisting of the equivalence relations.

Given an equivalence $\varepsilon$ on $X$, we denote by $<\varepsilon>$ the smallest coarse structure on $X$
such that $\varepsilon\in <\varepsilon>$  and
say that $<\varepsilon>$ is {\it generated} by $\varepsilon$.
To show that $<\varepsilon>$ is  cellular, we take the partition $\mathcal{P}$ into $\varepsilon$-classes and, for each finite  subset  $\mathfrak{F}$ of $\mathcal{P}$, denote by $\mathcal{P}_{\mathfrak{F}}$
the partition $\cup\mathfrak{F}$,  $P$, $P\in\mathcal{P}\setminus\mathfrak{F}$. We denote by $\varepsilon_{\mathfrak{F}}$ the equivalence defined by the partition $\mathcal{P}_{\mathfrak{F}}$ and observe that
$\{\varepsilon_{\mathfrak{F}}: \mathfrak{F}$
 is a finite subset of $\mathcal{P}\}$  is a base for $<\varepsilon>$.

By [9, Theorem 3.11],  a metrizable coarse structure $\mathcal{E}$ on $X$ is cellular if and  only if  $\mathcal{E}=\mathcal{E}_{d}$  for some
 ultrametric $d$ on $X$.  By  [9, Theorem 3.1.3],   a  coarse
 structure $\mathcal{E}$  on $X$ is cellular  if and only if $asdim (X, \mathcal{E})=0$.
 For definition of the asymptotic dimension see  [9, p.45] or  [11, p. 129].

If $\mathcal{E}_{1}, \mathcal{E}_{2}$  are cellular  coarse structure on $X$  then
$\mathcal{E}_{1}\wedge \mathcal{E}_{2}$  is cellular,  but the set  of  all cellular coarse
 structure on $X$ is not closed under the operation $\vee$.
 Moreover,  at the  moment we may even conjecture that every coarse structure $\mathcal{E}$ on $X$
  can be represented as $\mathcal{E}=\mathcal{E}_{1}\vee \mathcal{E}_{2}$
    for some cellular coarse structures on $X$. After two examples,
     we put this insolent conjecture in restricted form  of open questions.

We say that a coarse structure $\mathcal{E}$  on $X$  is $\vee$-{\it decomposable} if there exist two cellular coarse structures $\mathcal{E}_{0}, \mathcal{E}_{1} $
 on $X$ such that $\mathcal{E}=\mathcal{E}_{0}\vee \mathcal{E}_{1}$.
\vspace{4 mm}

{\bf Example 1.}
For  $n\in \mathbb{N}$, we endow $\mathbb{Z}^{n}$ with the metric
$ d(x,y)= \mid x_{1} -y_{1}\mid +\ldots   +$ $\mid x_{n} - y_{n}\mid $,
 put  $\mathcal{E}=\mathcal{E}_{d}$ and show
that $\mathcal{E}$ is $\vee$-decomposable.

We put  $K=\{x\in \mathbb{Z}^{n}:  x_{i}\in \{0,1\}\}$ and consider the partitions
$\mathcal{P}_{0}= \{ K + a: a\in 2 \mathbb{Z}^{n}\}$, $\mathcal{P}_{1}=\{(1, 1, \ldots , 1)+ P:  P\in\mathcal{P}_{0}\}$.
Let $\varepsilon_{0}, \varepsilon_{1}$ be the equivalences defined by $\mathcal{P}_{0}$,  $\mathcal{P}_{1}$
and $\mathcal{E}_{0}$,  $\mathcal{E}_{1}$
be the coarse structures generated by $\varepsilon_{0}$,  $\varepsilon_{1}$. Clearly,
$\mathcal{E}_{0}\subset \mathcal{E}$,  $\mathcal{E}_{1}\subset\mathcal{E}$ so $\mathcal{E}_{0}\vee \mathcal{E}_{1}
\subseteq\mathcal{E}$.
If   $x, y\in\mathbb{Z}^{n}$ and $d(x,y)=1$
then there exists $P\in \mathcal{P}_{0}\cup \mathcal{P}_{1}$
such that $\{x,y\}\subset P$.
It follows that $ \mathcal{E}\subseteq \mathcal{E}_{0}\vee\mathcal{E}_{1}$.

\vspace{4 mm}

Let $\Gamma$ be a connected graph with the set of vertices $\mathcal{V}$.
 We endow $\mathcal{V}$  with the path metric $d$  and denote
 by $\mathcal{E} _{\Gamma}$ the coarse structure $\mathcal{E} _{d}$ on $\mathcal{V}$.
 Given a coarse structure  $\mathcal{E}$ on a set $X$,  how one can detect if there exists a graph $\Gamma$
 with the set of vertices $X$ such that $\mathcal{E}=\mathcal{E}_{\Gamma}$?
 The answer to this question gives Theorem 5.1.1 from [9].
\vskip 7pt

\vspace{4 mm}

 {\bf Example 2.} Let $\mathcal{T}$ be a tree. Then the coarse structure $\mathcal{E}_{\mathcal{T}}$ is
 $\vee$-decomposable, see Proposition 9 for  more general result.

\vspace{4 mm}

{\bf Question 1.} {\it Is every uniformly locally finite coarse structure  $\vee$-decomposable?}

\vspace{4 mm}

{\bf Question 2.} {\it Is the coarse structure $\mathcal{E}_{\Gamma}$  $\vee$-decomposable
for every connected graph  $\Gamma$?}
\vspace{4 mm}

We note  that the finitary coarse structures of finitely generated  groups lie in the  intersection of Questions 1 and 2.

Let $(\Gamma_{n})_{n\in\omega}$ be a sequence of finite connected  graphs
 with pairwise disjoint sets of vertices  $\{\mathcal{V}_{n}: n\in \omega\}$. For
  each $n\in\omega$, we pick $\upsilon _{n} \in \mathcal{V}_{n}$  and join $\Gamma_{n}$ and $\Gamma_{n+1}$  by
  the  edge $\{\upsilon_{n}, \upsilon_{n+1}\}$.  We believe that the resulting graph $\Gamma$ is a good candidate for  counterexamples to both questions if $(\Gamma_{n}) _{n\in\omega}$ is a family  of expanders. What is an expander, see Wikipedia.

\vskip 7pt

{\bf Proposition 9. } {\it Let $\Gamma$  be a graph with the set of  vertices $\mathcal{V}$, $v_{0}\in\mathcal{V}$,  $S(v_{0}, n)=\{v\in\mathcal{V}: d(v_{0}, v)=n\}$.
We denote by $\Gamma_{n}$
the subgraph of $\Gamma$ with the set of vertices $S(v_{0}, n)\cup S(v_{0}, n+1)$
and assume that there exists $k\in\mathbb{N}$
 such that, for  every $n\in\omega$,
 each connected component of $\Gamma_{n}$ is contained in some ball of radius $k$ in $\Gamma$.

Then the coarse structure $\mathcal{E}_{\Gamma}$ is $\vee$-decomposable.
\vskip 7pt

Proof.}
We denote by $\mathcal{P}_{0}$ the partition of $\mathcal{V}$ into
 connected components of $\Gamma_{2n}$, $n\in\omega$, and by $\mathcal{P}_{1}$
  the partition of $\mathcal{V}\setminus\{v_{0}\}$
  into connected components of
  $\Gamma_{2n+1}$,  $n\in\omega$.
Then we take the equivalences $\varepsilon_{0}, \varepsilon_{1}$
 on $\mathcal{V}$   defined by the partitions $\mathcal{P}_{0}$ and $\mathcal{P}_{1}$, $\{v_{0}\}$,
   and denote by  $\mathcal{E}_{0}$, $\mathcal{E}_{1}$
   the coarse structures on $\mathcal{V}$ generated by $\mathcal{E}_{0}$ and  $\mathcal{E}_{1}$.
By the assumption, $\mathcal{E}_{0} \subseteq \mathcal{E}_{\Gamma}$,  $\mathcal{E}_{1} \subseteq \mathcal{E}_{\Gamma}$
so $\mathcal{E}_{0} \vee \mathcal{E}_{1} \subseteq \mathcal{E}_{1}$.
Since every edge of $\Gamma$ is an edge of some connected component of $\Gamma_{m}$,  we have
$\mathcal{E}_{\Gamma} \subseteq \mathcal{E}_{0} \vee \mathcal{E}_{1}$.  $ \ \  \  \  \Box$
\vskip 7pt

We note that beside trees, Proposition 9 demonstrates $\vee$-decomposability of $\mathcal{E}_{\Gamma}$
provided that there exists $k\in\mathbb{N}$  such that $|S(v_{0},n)|\leq k$
 for each $n\in \omega$.
\vskip 7pt

{\bf Proposition 10. } {\it
   Assume that there exists a partition $\mathcal{P}$ of the set $\mathcal{V}$ of vertices of a graph $\Gamma$
    and $r\in \mathbb{N}$ such that, for every $P\in \mathcal{P}$ the following conditions hold:
\vskip 4pt

$(i) \  \  $  $  \  \  P \subseteq B(v,r)$  for  some $v\in \mathcal{V}$;

$(ii) \  \  $   $ \  \  \mid\{P^{\prime}\in\mathcal{P}:  P\neq P^{\prime},  \  \
P^{\prime}\cap B(P,1)\neq\emptyset \}\mid\leq \mid P\mid$.
\vskip 4pt

Then the coarse structure  $\mathcal{E}_{\Gamma}$
is  $\vee$-decomposable.

\vskip 7pt

Proof.}
We denote $\mathcal{N}(P)= \{ P^{\prime}\in \mathcal{P}: P\neq P^{\prime}, P^{\prime}\cap B(P,1)\neq\emptyset\}$
and suppose that we have chosen a family $\mathfrak{F}$ of two-element pairwise disjoint subsets of $\mathcal{V}$  such that
\vskip 4pt

$(*) \  \  $ for any $P\in\mathcal{P}$  and $Q\in\mathcal{N}(P)$,  there exists
$F\in\mathfrak{F}$ such that $\mid P\cap F\mid = \mid Q \cap F\mid=1$,
 and for any $F\in\mathfrak{F}$,
  there exist $P\in\mathcal{P}$,  $Q\in\mathcal{N}(P)$  such that
  $\mid P\cap F\mid = \mid Q \cap F\mid=1$.
\vskip 5pt

We take the equivalences $\varepsilon, \gamma$ on  $\mathcal{V}$ defined
by the partitions $\mathcal{P}$ and $\mathfrak{F}$, $\{v\}$,  $v\in\mathcal{V}\setminus \cup\mathfrak{F}$,
  and denote by $\mathcal{E}_{0}$
  and $\mathcal{E}_{1}$  the coarse structures on   $\mathcal{V}$  generated by  $\varepsilon$ and $\gamma$.
We observe that $\mathcal{E}_{0}$, $\mathcal{E}_{1}$     are cellular and, by $(i)$, $(\ast)$,
$\mathcal{E}_{0}\subseteq \mathcal{E}_{\Gamma}$,   $\mathcal{E}_{1}\subseteq \mathcal{E}_{\Gamma}$  so
$\mathcal{E}_{0}\vee\mathcal{E}_{1}\subseteq \mathcal{E}_{\Gamma}$.

On the other hand,  let $\{x , y\}$  be an edge of $\Gamma$.
If $\{x,y\}\subseteq P$ for some $P\in\mathcal{P}$, then
$(x,y)\in\varepsilon$. Otherwise, there exist $P\in\mathcal{P}$
 and $Q\in\mathcal{N}(P)$ such that
 $x\in P$, $y\in Q$ and, by $(\ast)$, $(x,y)\in\varepsilon\circ\gamma\circ\varepsilon$.
Hence,  $\mathcal{E}_{\Gamma}\subseteq \mathcal{E}_{0}\vee\mathcal{E}_{1}$.

To define $\mathfrak{F}$, we enumerate $\mathcal{P}=\{P_{\alpha}: \alpha< \kappa\}$
 and choose an injective mapping $f_{\alpha}: \mathcal{N}(P_{\alpha})\longrightarrow P_{\alpha}$.
  We denote
  $X_{\alpha}=\{\beta<\alpha: P_{\beta}\in\mathcal{N}(P_{\alpha})\}$ and put
  $\mathfrak{F}_{\alpha}=\{\{f_{\beta}(P_{\alpha}), f_{\alpha}(P_{\beta})\}: \beta\in X_{\alpha}\}$.
   Then the desired $\mathfrak{F}$ is $\bigcup_{\alpha< \kappa} \mathfrak{F}_{\alpha}$. $ \ \ \  \Box$
 \vskip 7pt

{\bf Proposition 11. } {\it
Let $\Gamma$ be a graph  with the set of vertices $\mathcal{V}$.
  Assume that  there exists $r\in \mathbb{N}$ such that $|B(v,r)|= |\mathcal{V}|$
   for each $v\in \mathcal{V}$. Then the coarse structure $\mathcal{E}_{\Gamma}$ is $\vee$-decomposable.
\vskip 4pt

Proof.} We use the Zorn's lemma to choose a subset $X \subset \mathcal{V}$
 such that $B(x,r)\cap B(x^{\prime},r)=\emptyset$  for all distinct
 $x, x^{\prime} \in X$    and, for every $v\in \mathcal{V}$,
 there exists $x\in X$  such that
 $B(v,r)\cap B(x,r)\neq\emptyset$.
Then we take an arbitrary partition $\mathcal{P}$ of $\mathcal{V}$ such that if  $x\in X$ and
$x\in P$, $P\in \mathcal{P}$ then
$P\subseteq B(x, 2r)$. Apply Proposition 10.
$ \ \ \  \Box$
 \vskip 7pt

We recall that a metric space $(X, d)$ is {\it geodesic} if,
for any $x,y\in X$, there is an isometric embedding $f: [0, d(x,y)]\longrightarrow X$
 such that $f(0)=x$, $f(d(x,y))=y$.
 \vskip 7pt

 {\bf Proposition 12. } {\it
Let $(X,d)$ be a geodesic metric space such that   $|B_{d}(x,1)|=X$   for each $x\in X$. Then the
coarse structure $\mathcal{E}_{d}$ is $\vee$-decomposable.

\vskip 4pt

Proof.}
We consider a graph $\Gamma$  with the set of vertices $X$ and  the set of edges
$\{\{x,y\}: 0< d(x,y)\leq 1\}$.
Then the coarse structures  $\mathcal{E}_{d}$  and $\mathcal{E}_{\Gamma}$  coincide.
Apply Proposition 11.$ \ \ \  \Box$
 \vskip 7pt

In all above propositions, the coarse structure $\mathcal{E}_{0}, \mathcal{E}_{1}$
 witnessing the decomposition $\mathcal{E}=\mathcal{E}_{0}\vee\mathcal{E}_{1}$  are generated
  by some equivalences $\varepsilon_{0}$, $\varepsilon_{1}$.
If in Proposition 12 $(X, d)$ is separable then $\mathcal{E}_{0}, \  \mathcal{E}_{1}$ can be chosen to be metrisable.

\vspace{4 mm}

Let $(X, \mathcal{E})$ be   a coarse space. Each non-empty subset  $Y$  of $X$  has the natural coarse structure
$\mathcal{E}| _{Y}= \{\varepsilon\cap(Y\times  Y): \varepsilon\in \mathcal{E}\}$
 which is called a {\it substructure} of  $\mathcal{E}$.

\vspace{4 mm}

A subset $Y$  of $X$  is called {\it large } in  $(X, \mathcal{E})$ if there exists  $\varepsilon\in \mathcal{E}$  such that $X=B(Y, \mathcal{E})$.

\vspace{4 mm}

{\bf Proposition 13.} {\it
Let $\mathcal{E}$  be a coarse structure on a set $X$  and let $Y$ be a large subset of $X$.
If $\mathcal{E}| _{Y}$ is decomposable then $\mathcal{E}$ is  $\vee$-decomposable.

\vspace{4 mm}

Proof.} We take $\varepsilon_{0}\in \mathcal{E}$ such that $\varepsilon_{0}=\varepsilon_{0}^{-1}$,  $X=B(Y, \varepsilon_{0})$ and
choose a mapping $f: X \longrightarrow Y$  such that $(x, f(x))\in \varepsilon_{0} $  and
$f(y)=y$ for every $y\in Y$. For each $\gamma\subseteq Y\times Y$, we put
$f^{-1}(\gamma)=$ $\{(x, x^{\prime})\in X\times X:  (f(x), f(x^{\prime}) \in\gamma)\}$ and   note that
\vspace{4 mm}

$(\ast)$   $\mathcal{U}$ is a base of $\mathcal{E}| _{Y}$ if and only if
$f^{-1}(\mathcal{U})=$ $\{f^{-1}(\gamma): \gamma\in\mathcal{U}\}$ is a base of $\mathcal{E}$.
\vspace{4 mm}

We assume that  $\mathcal{E}| _{Y}= \mathcal{E}_{0}^{\prime} \vee \mathcal{E}_{1}^{\prime} $
with cellular $\mathcal{E}_{0}^{\prime}, \mathcal{E}_{1}^{\prime} $
and put  $\mathcal{E}_{0}= f^{-1}( \mathcal{E}_{0}^{\prime})$,   $\mathcal{E}_{1}= f^{-1}( \mathcal{E}_{1}^{\prime})$
 and note that $\mathcal{E}_{0}, \mathcal{E}_{1} $
   are cellular. If  $\gamma$, $\gamma^{\prime}\in Y\times Y$  then
   $f^{-1}(\gamma\circ\gamma^{\prime})=$  $f^{-1}(\gamma)\circ f^{-1}(\gamma^{\prime})$
    so $f^{-1}(\mathcal{E}_{0}^{\prime} \vee \mathcal{E}_{1}^{\prime} )=  $ $\mathcal{E}_{0} \vee \mathcal{E}_{1}$
    and,  applying $(\ast)$, we get  $\mathcal{E}= \mathcal{E}_{0} \vee \mathcal{E}_{1}$.  $  \  \   \  \  \Box$

    \vspace{6 mm}

{\bf Question 3.} {\it Is every substructure of $ \ \ \vee$-decomposable coarse structure $\vee$-decomposable?}
\vspace{4 mm}

In view of Proposition 12, in order to  answer Question 2 in affirmative, it suffices to get the positive answer to Question 3 for substructures  $\mathcal{E}| _{Y}$ where $Y$ is a large subset of $X$. Indeed, replacing each edge $\{x,y\}$  of $\Gamma$ by the isometric copy of the unit interval, we get the geodesic metric space $(X,d)$ in which the set $\mathcal{V}$  of vertices of $\Gamma$ is large.

\vspace{4 mm}

Each ideal $\mathcal{I}$ in the  Boolean algebra   $\mathcal{P}_{X}$ containing $\mathfrak{F}_{X}$
defines  a cellular coarse structure $\mathcal{E}_{\mathcal{I}}$  with the base
$\{\varepsilon_{A} : A\in \mathcal{I}\}$  where $\varepsilon_{A}$   is the partition $A$, $\{x\}$,   $x\in X \setminus A$.
Thus,  $B(x, \varepsilon_{A})=A$ if $x\in A$  and $B(x, \varepsilon_{A})= \{x\}$  if $x\in X \setminus  A$.  If $A$ is
an arbitrary subset of  $X$ then the set   $\mathcal{I}_{A}=\{A\cup  F: F\in \mathfrak{F}_{X}\}$  is the smallest   ideal such that  $A\in \mathcal{I}_{A}$  and $\mathfrak{F}_{X}\subseteq  \mathcal{I}_{A}$.

Given an arbitrary coarse structure $\mathcal{E}$ on $X$,
$\mathcal{E}\neq {\bf 1}_{X}$,  the family $\mathcal{I}$ of all bounded subset of  $X$
is   an ideal in  $\mathcal{P}_{X}$. We say that $\mathcal{E}_{\mathcal{I}}$ is the {\it companion} of $\mathcal{E}$.

An ideal $\mathcal{I}$ in $\mathcal{P}_{X}$ is maximal if and only if
$\{X \backslash A: A \in \mathcal{I}\}$   is a free ultrafilters, so
 there  are $2^{2^{\mid X\mid}}$  coarse structures on $X$ defined by  the maximal ideals, in particular,
 $| \mathcal{L}_{X} |= 2^{2^{\mid X\mid}}$.

\vspace{4 mm}

{\bf Proposition 14.} {\it
Let $\mathcal{E}$  be a cellular coarse structure on $X$ such that
$\mathcal{E}\neq {\bf 0}_{X}$
 and $\mathcal{E}\neq {\bf 1}_{X}$.
Then there  exists  a cellular coarse structure $\mathcal{E}^{\prime}$ on $X$  such
 that  $\mathcal{E}\vee\mathcal{E}^{\prime}= {\bf 1}_{X} $  and $\mathcal{E}^{\prime}\neq{\bf 1}_{X}$.
 \vspace{4 mm}

Proof.} By the assumption, there exists an equivalence
$\varepsilon\in \mathcal{E}$  such that either some $\varepsilon$-class $P$ is  infinite
 and  $X\setminus P$ is infinite,  or there are infinitely many
 $\varepsilon$-classes containing at least two elements.

In both cases, we choose a subset $A$ of $X$ such that $| A\cap B| =1 $   for each  $\varepsilon$-class  $B$,
put  $\mathcal{E}^{\prime}=\mathcal{E} _{\mathcal{I}_{A}}$
and observe that $\mathcal{E}\vee\mathcal{E}^{\prime}= {\bf 1}_{X} $,  $\mathcal{E}^{\prime}\neq {\bf 1}_{X} $.
$ \  \Box$

\vspace{4 mm}

{\bf Proposition 15.} {\it
Let $\mathcal{I}$ be a maximal ideal in $\mathcal{P}_{X}$ and $\mathcal{E}=\mathcal{E}_{\mathcal{I}}$.
If  $\mathcal{E}^{\prime}$ is a coarse structure on  $X$ such that $\mathcal{E}\wedge\mathcal{E}^{\prime}= {\bf 0}_{X} $ then $\mathcal{E}^{\prime}= {\bf 0}_{X} $.

\vspace{4 mm}

Proof.} We assume that  some infinite subset of $X$ is bounded in $\vee$.
Since the set
$\varphi=\{X\backslash A: A \in\mathcal{I}\}$  is a free
ultrafilter, there is an infinite subset $A$ of $X$ bounded in
$\mathcal{E}^{\prime}$ such that $A\in \mathcal{I}$.
We take $\gamma\in\mathcal{E}^{\prime}$
 such that $A\subseteq B(x,\gamma)$  for each $x\in A$. For every
 finite subset $K$ of $X$  and $x\in A\setminus K$, the set
  $A \cap B(x, \gamma)$       is infinite.
Hence, $\mathcal{E}\wedge\mathcal{E}^{\prime}\neq {\bf 0}_{X}$. Thus,
 every bounded in $\mathcal{E}^{\prime}$ subset of $X$ is finite.

Given an arbitrary $\gamma\in\mathcal{E}^{\prime}$, we denote $C=\{x\in X: \mid  B(x,\alpha)\mid > 1\}$
If $C$ is finite then
$\mathcal{E}^{\prime}= {\bf 0}_{X}$.
Otherwise $C$ is unbounded.
Since $\varphi$ is a free ultrafilter, there exists an infinite subset $C^{\prime}$  of $C$ such
that  $B(C^{\prime},\gamma ) \in  \mathcal{I}$.
By above paragraph,  $\mathcal{E}\wedge\mathcal{E}^{\prime}\neq {\bf 0}_{X}$.
$ \  \Box$

\vspace{4 mm}

Let $\mathcal{E}$ be a coarse structure on $X$.
We say that a coarse structure $\mathcal{E}^{\prime}$ is a {\it complement} to $\mathcal{E}$ if
$\mathcal{E}\wedge\mathcal{E}^{\prime}= {\bf 0}_{X}$, $\mathcal{E}\wedge\mathcal{E}^{\prime}= {\bf 1}_{X}$
If $\mathcal{E}$ has a complement, we say that $\mathcal{E}$ is complementable.
Clearly,  ${\bf 0}_{X}$ and ${\bf 1}_{X}$  are   {\it complementable}.
\vspace{4 mm}

 {\bf Example 3.} We take an infinite subset $A$ such that the set $B= X\setminus A$  is infinite, and denote
 $\mathcal{I}_{A}= \{A\cup F: F\in \mathfrak{F}_{X}\}$,  $\mathcal{I}_{B}= \{B\cup F: F\in \mathfrak{F}_{X}\}$.
Then $\mathcal{E}_{\mathcal{I}_{A}}\wedge   \mathcal{E}_{\mathcal{I}_{B}} = {\bf 0}_{X}$, $\mathcal{E}_{\mathcal{I}_{A}}\vee   \mathcal{E}_{\mathcal{I}_{B}} = {\bf 1}_{X}$.

\vspace{4 mm}

{\bf Question 4.} {\it How  one can detect  whether  a coarse structure on $X$ is complementable?}


\section{ Comments}

    1.	By Zorn's  lemma, every  unbounded coarse structure
on a set $X$ is contained is some  maximal unbounded      coarse structure.  For criterion of maximality   and some    properties of  maximal coarse structures see [9, Chapter 10].
We   mention only two facts. If $\mathcal{I}$  is a maximal ideal
on $X$ then $\mathcal{E} _{\mathcal{I}}$ is maximal. If $X$ is countable then the unit in the lattice of all uniformly finite coarse structure on $X$  (see Proposition 4) is maximal.
\vspace{4 mm}

2. We recall that an element $a$  of a lattice   $\mathcal{L}$ with
 ${\bf 0}$  and ${\bf 1} $ is an  {\it atom} (a {\it coatom}) if $a\neq {\bf 0}$ ($a\neq{\bf 1}$ ) and $b< a$ ($b>a$) implies $b= {\bf 0}$ ($b={\bf 1}$ ).

For every set $X$, the lattice  $\mathcal{L}_{X}$ has $2^{2 ^{|X|}}$ coatoms,
 see above paragraph.  We show that $\mathcal{L}_{X}$  has no atoms.
  Assume the contrary and let $\mathcal{E}$  be an atom in $\mathcal{L}_{X}$.

We suppose that the ideal $\mathcal{I}$  of all bounded subsets of $X$
has an infinite member $A$, partition $A$  into infinite subsets $B, C$ and denote by
$\mathcal{I}_{A}$ the smallest ideal in $\mathcal{P}_{X}$
containing $A$  and $\mathfrak{F}_{X}$. Then $\mathcal{E}_{\mathcal{I}_{A}} \subset  \mathcal{E}_{\mathcal{I}}\subseteq \mathcal{E}$, $\mathcal{E}_{\mathcal{I}_{A}} \neq {\bf 0}_{X}$.
It follows that $\mathcal{I}=\mathfrak{F}_{X}$.

Since $\mathcal{E}\neq {\bf 0}_{X}$,  we can choose $\varepsilon\in\mathcal{E}$
and a sequence $(a_{n})_{n\in\omega}$
in $X$ such that $|B(a_{n}, \varepsilon)|>1$, $B(a_{n}, \varepsilon) \cap B(a_{m}, \varepsilon)=\emptyset$
 for all distinct $n, m$.  For each $n\in \omega$,  we take $b_{2n}\in B(a_{2n}, \varepsilon)$,
 $b_{2n}\neq a_{2n}$ and denote $\gamma= \{(a_{2n}, b_{2n}), (b_{2n}, a_{2n}): n\in\omega\}\cup\bigtriangleup_{X}$.
We denote by $\mathcal{E}^{\prime}$  the smallest coarse structure on $X$ containing $\gamma$:
By the choice of $\gamma$ and above paragraph $(\mathcal{I}=\mathfrak{F}_{X})$,  ${\bf 0}_{X}\subset\mathcal{E}^{\prime}$ and, for every $\delta\in\mathcal{E}^{\prime}$,
 there is $m\in \omega$  such that $|B(a_{2m+1}, \delta)|=1$.
Hence, ${\bf 0}_{X}\subset\mathcal{E}^{\prime}\subset \mathcal{E}$  and we run in a contradiction.

\vspace{4 mm}

If $\mathcal{E}$  is a metrizable coarse structure on $X$ and $\mathrm{E}\neq {\bf 0}_{X} $  then Proposition 1 and these arguments give us metrizable $\mathcal{E}^{\prime}$ such that  ${\bf 0}_{X}\subset  \mathcal{E}^{\prime} \subset  \mathcal{E}$.

\vspace{4 mm}
3.	For a group $G$, the family $LI_{G}$ of all group ideals in
 $\mathcal{P}_{G}$ has the natural lattice structure.
 The lattices $LI_{G}$ were introduced and studied in [3],
  see also  [9, Chapter 6]. We note that the lattice $LI_{G}$
   is isomorphic to the lattice of all group coarse
   structure on $G$ (or group balleans in terminology of [9]).
Every group coarse structure $\mathcal{E}$ on  $G$ has the base
 $\{(x,y): x^{-1} y  \in A\}$, $A\in\mathcal{I}$  where $\mathcal{I}$ is a group ideal
  in $\mathcal{P}_G$.

For every Abelian group $G$, the lattice $LI_{G}$  is
modular, but $LI_{F}$  is not modular for a free
 non-Abelian group $F$.
It would be interesting to find a lattice-theoretical property valid for each $LI_{G}$   and, more generally, for every lattice $\mathcal{L}_{X}$.
\vspace{4 mm}

4.  By    analogy with topological groups, a group  coarse structure $\mathcal{E}$ on a group $G$ is called maximal $\mathcal{E}$ is a maximal unbounded coarse structure on $G$. These structures were introduced and studied in [6], see  also [9, Section 10.3]. If a group $G$ admits a maximal group coarse structure $\mathcal{E}$  then  the set $\{y^{2}: g\in G\}$  is bounded in $\mathcal{E}$. It follows that, say, $\mathbb{Z}$ does not admit such a structure.  Under $CH$,  the maximal group  coarse structure on the countable Boolean group was constructed  in [6], see also [9, Example 10.3.2].  As to our knowledge, the question whether or not, a  maximal group coarse structure can be constructed in ZFC with  no additional  set-theoretical assumptions remains open.
\vspace{4 mm}

5. Let $\Gamma$ be a connected graph  with the set of vertices $\mathcal{V}$.
 We replace each edge $\{u,v\}$ of $\Gamma$ to  three edges
  $\{u, x _{\{u,v\}}\}, \  \  \{x _{\{u,v\}},  y _{\{u,v\}}  \}, \  \  \{y _{\{u,v\}},  y\}$
  and get the subdivision $\Gamma^{\prime}$
    of  $\Gamma$. We denote by $d$
    the pass metric on the set $\mathcal{V}^{\prime}$  of vertices of  $\Gamma^{\prime}$.
The partition  $\{B_{d}(v,1):  v\in \mathcal{V}\}$
of $\mathcal{V}^{\prime}$  satisfies Proposition 10 so
$\mathcal{E}_{\Gamma^{\prime}}$ is $\vee$-decomposable.

We note that  $\mathcal{E}_{\Gamma}$ is uniformly locally finite if
and only if there exists $m\in \mathbb{N}$  such that
$\lambda(v)\leq m$ for each $v\in \mathcal{V}$, where  $\lambda(v) $ is the local degree  of $v$.
By above paragraph, if $\mathcal{E}_{\Gamma}$  is uniformly locally finite then there exist a set
$\mathcal{V}^{\prime}$,  $\mathcal{V} \subset\mathcal{V}^{\prime}$
and  a uniformly locally finite structure $\mathcal{E}^{\prime}$ on $\mathcal{V}^{\prime}$  such that
 $\mathcal{E}^{\prime}$  is   $\vee$-decomposable, $\mathcal{V}$ is a large subset  of  $(\mathcal{V}^{\prime}, \mathcal{E}^{\prime})$  and $\mathcal{E}^{\prime}|_{\mathcal{V}} = \mathcal{E}_{\Gamma} $.
 \vspace{4 mm}

6. {\it Can every metric uniformity on a set $X$ be represented as supremum of two ultrametric uniformities?}
As to our knowledge, this question was nowhere asked.

We show that, for every set $X$, the metric uniformity $\mathcal{M}$ of $\ell_{\infty}(X)$ is supremum of two ultrametric uniformities.  We denote
$P=\{f\in \ell_{\infty}(X): 0\leq f<1\}, $  $ R= \{f\in \ell_{\infty}(X): -\frac{1}{2}\leq f<\frac{1}{2}\}$,
 transfer $P$ and $R$ by all integer valued functions from $\ell_{\infty}(X)$
 and get the partitions  $\mathcal{P}$  and $\mathcal{R}$  of $ \ell_{\infty}(X)$.
For each $n\in \omega$, we denote
$\mathcal{P}_{n}= \{(\frac{1}{2})^{n} P^{\prime} : P^{\prime}\in \mathcal{P}\}$,
$\mathcal{R}_{n}= \{(\frac{1}{2})^{n} R^{\prime} : R^{\prime}\in \mathcal{R}\}$. Then
$\{\mathcal{P}_{ n}: n\in\omega\}$
 and $\{\mathcal{R}_{ n}: n\in\omega\}$
 define   ultrametric uniformities $\widetilde{\mathcal{P}}$
  and $\widetilde{\mathcal{R}}$
  such that $\mathcal{M}=\widetilde{\mathcal{P}} \ \widetilde{\mathcal{R}}$, in particular,
  $\mathcal{M}=\sup\{\widetilde{\mathcal{P}}, \ \widetilde{\mathcal{R}}\}$.

If  we use the families of partitions
$\{\mathcal{P}_{n}^{\prime}: n\in\omega\}$, $\{\mathcal{R}_{n}^{\prime}: n\in\omega\}$,
 where
 $\mathcal{P}_{n}^{\prime}= \{2^{n} P^{\prime}: P^{\prime}\in\mathcal{P}\}$,
 $\mathcal{R}_{n}^{\prime}= \{2^{n} R^{\prime}: P^{\prime}\in\mathcal{R}\}$,
 then get $\vee$-decomposition of the coarse structure of $ \ell_{\infty}(X)$
  in two cellular.

Every metric space $(X, d)$ can be realized isometrically as some subspace of $ \ell_{\infty}(X)$ but we do not know
 if decomposability  is inherited by subspaces of $ \ell_{\infty}(X)$  neither  in coarse nor in uniform cases.


\vskip 5pt

CONTACT INFORMATION

I.~Protasov: \\
Faculty of Computer Science and Cybernetics  \\
        Kyiv University  \\
         Academic Glushkov pr. 4d  \\
         03680 Kyiv, Ukraine \\ i.v.protasov@gmail.com

\medskip

K.~Protasova:\\
Faculty of Computer Science and Cybernetics \\
        Kyiv University  \\
         Academic Glushkov pr. 4d  \\
         03680 Kyiv, Ukraine \\ ksuha@freenet.com.ua

\end{document}